\PassOptionsToPackage{unicode}{hyperref}
\PassOptionsToPackage{hyphens}{url}
\documentclass[
]{article}
\usepackage{xcolor}
\usepackage[margin=1in]{geometry}
\usepackage{amsmath,amssymb}
\usepackage{iftex}
\ifPDFTeX
  \usepackage[T1]{fontenc}
  \usepackage[utf8]{inputenc}
  \usepackage{textcomp} % provide euro and other symbols
\else % if luatex or xetex
  \usepackage{unicode-math} % this also loads fontspec
  \defaultfontfeatures{Scale=MatchLowercase}
  \defaultfontfeatures[\rmfamily]{Ligatures=TeX,Scale=1}
\fi
\usepackage{lmodern}
\ifPDFTeX\else
\fi
\IfFileExists{upquote.sty}{\usepackage{upquote}}{}
\IfFileExists{microtype.sty}{% use microtype if available
  \usepackage[]{microtype}
  \UseMicrotypeSet[protrusion]{basicmath} % disable protrusion for tt fonts
}{}
\makeatletter
\@ifundefined{KOMAClassName}{% if non-KOMA class
  \IfFileExists{parskip.sty}{%
    \usepackage{parskip}
  }{% else
    \setlength{\parindent}{0pt}
    \setlength{\parskip}{6pt plus 2pt minus 1pt}}
}{% if KOMA class
  \KOMAoptions{parskip=half}}
\makeatother
\providecommand{\tightlist}{%
  \setlength{\itemsep}{0pt}\setlength{\parskip}{0pt}}
\usepackage{bookmark}
\IfFileExists{xurl.sty}{\usepackage{xurl}}{} % add URL line breaks if available
\hypersetup{
  pdftitle={The Sampling Distribution of the Log-Euclidean Distance Between Sample Correlation Matrices},
  hidelinks,
  pdfcreator={LaTeX via pandoc}}

\title{The Sampling Distribution of the Log-Euclidean Distance Between
Sample Correlation Matrices}
\author{Argyn Kuketayev}
\date{\today}

\begin{document}
\maketitle

\section{The Sampling Distribution of the Log-Euclidean Distance Between
Sample Correlation
Matrices}\label{the-sampling-distribution-of-the-log-euclidean-distance-between-sample-correlation-matrices}

\subsection{Abstract}\label{abstract}

Comparing correlation matrices across time or stress scenarios is
critical in quantitative finance and multivariate statistics, yet sample
estimation noise often obscures whether an observed distance reflects a
true structural shift. We derive the asymptotic sampling distribution of
the intrinsic off-log (log-Euclidean) distance between two independently
estimated full-rank correlation matrices under the null hypothesis that
their population correlation matrices coincide. Under general sampling
with finite fourth moments, the scaled squared distance converges to a
weighted sum of independent \(\chi_1^2\) variables, with weights
determined by the asymptotic covariance of the Generalized Fisher
Transformation (GFT) coordinates. Under Gaussian sampling at
independence, this simplifies to a parameter-free \(4\chi_d^2\) law. To
calibrate tail probabilities, we provide closed-form cumulant generating
functions, Lugannani--Rice saddlepoint quantiles, and an explicit
Chernoff envelope requiring no root-finding. The first moment of the
limiting law establishes a simple rule of thumb for the baseline
expected distance under the null hypothesis
(\(\operatorname E[d_{\mathrm{LE}}] \lesssim 2\sqrt{d/n}\) near
independence), quantifying the average separation induced strictly by
estimation error. We establish plug-in consistency, present an explicit
Gaussian covariance factorization, compare the distance statistic with
coordinate Wald tests, and characterize its local power.

\begin{center}\rule{0.5\linewidth}{0.5pt}\end{center}

\subsection{1. Introduction}\label{introduction}

The statistical comparison of correlation structures is a foundational
problem across multivariate statistics, with applications in
neuroimaging and quantitative finance. In the latter setting, risk
managers often perturb correlations to study portfolio sensitivity, but
elementwise perturbations can violate positive definiteness. This has
motivated unconstrained parametrizations that map the manifold of
full-rank correlation matrices, \(\operatorname{Cor}_+(p)\), bijectively
into a Euclidean space \(\mathbb R^d\), where \(d=p(p-1)/2\).

A standard unconstrained parametrization is the Generalized Fisher
Transformation (GFT) of Archakov and Hansen (2021), \[
\gamma(C)=\operatorname{vecl}(\log C).
\] Thanwerdas and Pennec (2022) endowed \(\operatorname{Cor}_+(p)\) with
the flat off-log metric obtained by pulling back the Frobenius inner
product through the off-log coordinates, and Bisson and Pennec (2026)
developed the surrounding theory of log-Euclidean Lie groups. For the
standard log-Euclidean metric on \(\operatorname{SPD}(p)\), the
inclusion
\(\operatorname{Cor}_+(p)\hookrightarrow\operatorname{SPD}(p)\) is not
isometric; Bisson and Pennec (2026, Proposition 2.36) identify the extra
diagonal term. They also construct an adapted ambient log-Euclidean
metric for which the standard inclusion becomes an isometric, totally
geodesic embedding (Theorem 4.27). Intrinsically, the off-log geodesic
distance is simply Euclidean distance in GFT coordinates, up to the
duplication factor induced by symmetry.

Related geometric approaches have recently used geodesic distances to
construct correlation stress scenarios in financial risk management
(Chmielowski, 2026). That construction operates on a
constant-determinant submanifold of covariance matrices, along which
individual variances may vary, and grades scenarios by a plausibility
measure decreasing in Fisher--Rao distance. The off-log geometry instead
operates directly on the unit-diagonal correlation manifold, in the same
GFT coordinates whose sampling behavior is known.

While recent work such as Pereira, Mestre, and Gregoratti (2024) has
derived asymptotic distributions for log-Euclidean distances between
sample covariance matrices, the sampling law of the intrinsic distance
on the unit-diagonal correlation manifold has not been established. The
statistical literature supplies sampling theory at the level of
coordinates---such as the asymptotic covariance of \(\hat\gamma\) and
Wald tests built by whitening with it---but practitioners comparing
correlation structures generally report the intrinsic distance itself.

The structural link enabling this distance-level inference is that the
off-log distance is exactly \(\sqrt2\) times the Euclidean distance
between GFT coordinate vectors, because \(\psi\) acts as a global
isometry. Distance-level distribution theory is therefore an exact
consequence of coordinate-level theory. In this paper, we establish the
limiting distribution of the scaled squared distance under general
sampling with finite fourth moments, allowing for distinct fourth-moment
structures and unequal sample sizes across populations. We then develop
analytical tail-calibration tools, including the closed-form cumulant
generating function, a Lugannani--Rice approximation, and an explicit
Chernoff envelope. Under Gaussian sampling, we provide an operator
factorization of the GFT coordinate covariance and show that at
independence the limiting law reduces to \(4\chi_d^2\). Finally, we
analyze the local power of the distance statistic, showing that its
minimum detectable effect improves like \(d^{-1/4}\) for
coordinate-dense alternatives but requires sample sizes growing like
\(\sqrt d\) for systemic rank-one shifts.

Section 2 collects the geometric preliminaries. Section 3 states the
sampling setup and proves the main distance limit. Section 4 develops
the tail calibration of that limit. Section 5 gives the Gaussian
covariance factorization and pivotal specialization. Section 6 treats
unequal sample sizes, plug-in weights, the Wald comparison, local power,
and non-Gaussian sampling.

\begin{center}\rule{0.5\linewidth}{0.5pt}\end{center}

\subsection{2. Geometric Preliminaries and Log-Euclidean
Geodesics}\label{geometric-preliminaries-and-log-euclidean-geodesics}

Let \(\operatorname{SPD}(p)\) denote the cone of \(p\times p\) real
symmetric positive-definite matrices, and let \[
\operatorname{off}:\operatorname{Symm}(p)\to\mathbb R^d,\qquad d=\frac{p(p-1)}2,
\] extract the strictly lower-triangular entries of a symmetric matrix.
The forward GFT coordinate map is \[
\psi(C)=\operatorname{off}(\log C),\qquad C\in\operatorname{Cor}_+(p).
\] Archakov and Hansen (2021, Theorem 1) show that this map is a
bijection. Equivalently, if
\(S(\mathbf u)=\operatorname{off}^{-1}(\mathbf u)\) is the hollow
symmetric matrix with off-diagonal coordinates \(\mathbf u\), then there
exists a unique additive diagonal correction
\(\mathcal D(\mathbf u)\in\operatorname{Diag}(p)\) such that \[
\psi^{-1}(\mathbf u)=\exp\!\big(S(\mathbf u)+\mathcal D(\mathbf u)\big)
\] has unit diagonal. Bisson and Pennec (2026, Theorem 2.29) restate
this construction geometrically.

The off-log metric is the pullback of the Frobenius inner product
through the hollow-matrix off-log map (Thanwerdas and Pennec, 2022;
Bisson and Pennec, 2026, Definition 2.31). Therefore, by the general
log-Euclidean distance formula of Bisson and Pennec (2026, Proposition
2.18(c)), \[
d_{\mathrm{LE}}(C_1,C_2)=\sqrt{2}\,\|\psi(C_1)-\psi(C_2)\|_2.
\] The factor \(\sqrt2\) arises because the Frobenius norm of a
symmetric hollow matrix counts each off-diagonal entry twice, whereas
\(\operatorname{off}\) retains only one copy. Because
\(C\mapsto S(\psi(C))\) is a global isometry from
\(\operatorname{Cor}_+(p)\) onto the Euclidean vector space of hollow
symmetric matrices equipped with the Frobenius inner product, off-log
geodesics are unique, extend for all real parameter values, and reduce
to straight-line interpolation and extrapolation in GFT coordinates;
every finite coordinate displacement maps back to a full-rank
correlation matrix.

\emph{Remark (Universe-dependent coordinates).} The map \(\psi\) does
not commute with restriction to a principal submatrix: for an index set
\(A\subset\{1,\ldots,p\}\), \[
\log\!\left(C_{AA}\right)\;\ne\;\left(\log C\right)_{AA}
\] in general, because the matrix logarithm nets out indirect paths
running through all remaining variables. The GFT coordinates of a fixed
pair, and hence the off-log distance between two correlation structures
on \(A\), therefore depend on the ambient variable set in which they are
computed. The effect is not small: for a \(100\)-variable one-factor
correlation matrix with average correlation \(0.53\), the ten
coordinates of a fixed five-variable subset computed within the
five-variable universe exceed those read off the \(100\times100\)
logarithm by a factor of roughly five. The mechanism is the same
\(O(1/p)\) compression visible in the equicorrelation coordinate
\(\gamma=[\log(1+(p-1)\rho)-\log(1-\rho)]/p\). Consequently, distances
computed in different ambient dimensions are not comparable, and the
sampling theory below applies to two estimates formed over a common
variable set. For geometrically consistent comparison of correlation
matrices of genuinely different sizes, see the nested isometric
embeddings of Bisson and Pennec (2026, Corollary 4.36).

\begin{center}\rule{0.5\linewidth}{0.5pt}\end{center}

\subsection{3. Sampling Distribution of the Geodesic
Distance}\label{sampling-distribution-of-the-geodesic-distance}

\subsubsection{3.1 Sampling setup}\label{sampling-setup}

For \(r=1,2\), let \(X_{r1},\ldots,X_{rn_r}\in\mathbb R^p\) be iid
observations from a distribution \(P_r\) with mean \(\mu_r\),
nonsingular covariance matrix, and finite fourth moments. The two
samples are independent. Because sample correlations are invariant to
translations and positive diagonal rescalings, we work without loss of
generality with the population-standardized observations, so that the
population covariance matrix in each sample is a correlation matrix.
Under the null hypothesis, \[
H_0:\qquad \operatorname{Corr}(P_1)=\operatorname{Corr}(P_2)=C\in\operatorname{Cor}_+(p),
\] while the two populations may otherwise have different fourth-moment
structures.

Let \[
\hat\Sigma_r=\frac1{n_r}\sum_{t=1}^{n_r}(X_{rt}-\bar X_r)(X_{rt}-\bar X_r)^T,
\qquad C_r=\pi(\hat\Sigma_r),
\] where \[
\pi(\Sigma)=(I\circ\Sigma)^{-1/2}\Sigma(I\circ\Sigma)^{-1/2}.
\] Write \(\mathbf u=\psi(C_1)\), \(\mathbf v=\psi(C_2)\), and
\(\boldsymbol\mu_C=\psi(C)\). For the population-standardized
observations define \[
\mathbf\Sigma_X^{(r)}
:=\operatorname{Var}\!\left[\operatorname{vec}\!\left((X_{r1}-\mu_r)(X_{r1}-\mu_r)^T\right)\right].
\] Finite fourth moments imply the usual multivariate central limit
theorem \[
\sqrt{n_r}\left(\operatorname{vec}(\hat\Sigma_r)-\operatorname{vec}(C)\right)
\xrightarrow{d}
\mathcal N\!\left(0,\mathbf\Sigma_X^{(r)}\right).
\] Equality of the population correlation matrices does not require
\(\mathbf\Sigma_X^{(1)}=\mathbf\Sigma_X^{(2)}\).

\subsubsection{3.2 Lemma (Asymptotic normality of GFT
coordinates)}\label{lemma-asymptotic-normality-of-gft-coordinates}

\emph{(Archakov and Hansen, 2021, Section 3.4).} Under the setup above,
\[
\sqrt{n_r}\big(\psi(C_r)-\psi(C)\big)
\xrightarrow{d}
\mathcal N\!\left(0,\mathbf V_C^{(r)}\right),\qquad r=1,2,
\] independently, where \[
\mathbf V_C^{(r)}
=J_\psi(C)J_\pi(C)\mathbf\Sigma_X^{(r)}J_\pi(C)^T J_\psi(C)^T.
\] Here \(J_\pi(C)\) is the Jacobian of the covariance-to-correlation
map and \(J_\psi(C)\) is the Jacobian of the GFT map. Classical explicit
formulas for the intermediate asymptotic covariance of the sample
correlation matrix are given by Neudecker and Wesselman (1990) and
Browne and Shapiro (1986).

\emph{Proof.} Both \(\pi\) and \(\psi\) are continuously differentiable
in a neighborhood of \(C\); in particular, the matrix logarithm is
analytic on \(\operatorname{SPD}(p)\). The claim follows from the delta
method applied separately to the two independent samples. \(\square\)

\subsubsection{3.3 Theorem (Equal-sample-size distance
law)}\label{theorem-equal-sample-size-distance-law}

Assume \(n_1=n_2=n\). Under \(H_0\), \[
\boxed{
 n\,d_{\mathrm{LE}}^2(C_1,C_2)
 \xrightarrow{d}
 \sum_{i=1}^d 2\lambda_i\!\left(\mathbf V_C^{(1)}+\mathbf V_C^{(2)}\right)\chi^2_{1,i}
 }
\] where the \(\chi^2_{1,i}\) are independent central chi-squared
variables and \(\lambda_i(\cdot)\) denotes the \(i\)th eigenvalue.

\emph{Proof.} Let \(\mathbf z=\mathbf u-\mathbf v\). Independence and
Lemma 3.2 give \[
\sqrt n\,\mathbf z\xrightarrow{d}\mathbf y\sim
\mathcal N\!\left(0,\mathbf V_C^{(1)}+\mathbf V_C^{(2)}\right).
\] Since \(d_{\mathrm{LE}}^2(C_1,C_2)=2\|\mathbf z\|_2^2\), \[
n\,d_{\mathrm{LE}}^2(C_1,C_2)
=2\|\sqrt n\,\mathbf z\|_2^2
\xrightarrow{d}2\mathbf y^T\mathbf y.
\] Diagonalizing \(\mathbf V_C^{(1)}+\mathbf V_C^{(2)}=Q\Lambda Q^T\)
yields \[
2\mathbf y^T\mathbf y\sim
\sum_{i=1}^d2\lambda_i\!\left(\mathbf V_C^{(1)}+\mathbf V_C^{(2)}\right)\chi^2_{1,i}
\] (Mathai and Provost, 1992, Chapter 3). \(\square\)

\subsubsection{3.4 Corollary (Common coordinate
covariance)}\label{corollary-common-coordinate-covariance}

If \(\mathbf V_C^{(1)}=\mathbf V_C^{(2)}=: \mathbf V_C\)---in
particular, if the two populations have the same relevant fourth-moment
structure---then \[
\boxed{
 n\,d_{\mathrm{LE}}^2(C_1,C_2)
 \xrightarrow{d}
 \sum_{i=1}^d4\lambda_i(\mathbf V_C)\chi^2_{1,i}.
 }
\] This is the form used in the Gaussian specialization below.

\subsubsection{3.5 Remark (Mean distance and the baseline
expectation)}\label{remark-mean-distance-and-the-baseline-expectation}

Throughout the remainder of the paper, write \[
\mathbf S_C:=\mathbf V_C^{(1)}+\mathbf V_C^{(2)},
\] so that the limiting weights of Theorem 3.3 are
\(w_i=2\lambda_i(\mathbf S_C)\), and \(\mathbf S_C=2\mathbf V_C\) in the
common-covariance case of Corollary 3.4. The first moment of the
limiting law and the induced bound on the limiting distance are \[
\operatorname E[Y]=2\operatorname{tr}(\mathbf S_C),
\qquad
\operatorname E[\sqrt Y]\le \sqrt{2\operatorname{tr}(\mathbf S_C)},
\] the second by Jensen's inequality. Under common covariance these
reduce to \(\operatorname E[Y]=4\operatorname{tr}(\mathbf V_C)\) and
\(\operatorname E[\sqrt Y]\le 2\sqrt{\operatorname{tr}(\mathbf V_C)}\).

Because \(d_{\mathrm{LE}}\) is a metric, the expected distance between
two independently estimated correlation matrices is strictly positive
under \(H_0\): estimation noise alone separates \(C_1\) and \(C_2\) on
the manifold, providing a baseline expectation for the distance under
the null hypothesis. The bound above therefore quantifies the average
separation induced strictly by estimation error. Writing \(d_\infty\)
for the limiting distance variable, so that
\(n\,d_{\mathrm{LE}}^2\to d_\infty^2\), we have
\(\operatorname E[d_\infty]\le2\sqrt{\operatorname{tr}(\mathbf V_C)}\)
under common covariance, and hence the finite-sample approximation \[
\operatorname E\!\left[d_{\mathrm{LE}}(C_1,C_2)\right]\;\lesssim\;2\sqrt{\frac{\operatorname{tr}(\mathbf V_C)}{n}}.
\] Near Gaussian independence
\(\operatorname{tr}(\mathbf V_C)\approx d\) by Section 5.3, giving the
simple rule of thumb
\(\operatorname E[d_{\mathrm{LE}}]\lesssim2\sqrt{d/n}\), with \(n\)
replaced by \(n_{\mathrm{eff}}\) when sample sizes are unequal.

\begin{center}\rule{0.5\linewidth}{0.5pt}\end{center}

\subsection{4. Tail Calibration of the Limiting
Law}\label{tail-calibration-of-the-limiting-law}

Applications of the distance statistic turn on the upper tail rather
than the center of the null law. For a significance level
\(\alpha\in(0,1)\) we require the critical distance \(y_\alpha\) with
\(P(Y>y_\alpha)=\alpha\), typically at small \(\alpha\). This section
collects the analytical machinery needed to obtain it. All expressions
are stated for the general weights \(w_i=2\lambda_i(\mathbf S_C)\) of
Theorem 3.3, with the common-covariance specialization of Corollary 3.4
recorded alongside.

\subsubsection{4.1 Cumulant generating function and its
derivatives}\label{cumulant-generating-function-and-its-derivatives}

For \(Y=\sum_{i=1}^dw_i\chi^2_{1,i}\) the cumulant generating function
is available in closed form: \[
K_Y(t)=-\frac12\sum_{i=1}^d\log\!\left(1-4\lambda_i(\mathbf S_C)t\right),
\qquad
 t<\frac1{4\lambda_{\max}(\mathbf S_C)},
\] reducing under common covariance to
\(K_Y(t)=-\tfrac12\sum_i\log(1-8\lambda_i(\mathbf V_C)t)\) for
\(t<1/(8\lambda_{\max}(\mathbf V_C))\). Its first two derivatives are
explicit: \[
K_Y'(t)=\sum_{i=1}^d\frac{2\lambda_i(\mathbf S_C)}{1-4\lambda_i(\mathbf S_C)t},
\qquad
K_Y''(t)=\sum_{i=1}^d\frac{8\lambda_i(\mathbf S_C)^2}{\left(1-4\lambda_i(\mathbf S_C)t\right)^2},
\] with common-covariance forms
\(\sum_i4\lambda_i(\mathbf V_C)/(1-8\lambda_i(\mathbf V_C)t)\) and
\(\sum_i32\lambda_i(\mathbf V_C)^2/(1-8\lambda_i(\mathbf V_C)t)^2\).

Since \(K_Y''>0\) on the convergence interval whenever
\(\mathbf S_C\ne0\), the function \(K_Y'\) is strictly increasing there,
with \(K_Y'(0)=\operatorname E[Y]\) and \(K_Y'(t)\to\infty\) as
\(t\uparrow1/(4\lambda_{\max}(\mathbf S_C))\). Hence for every
\(y>\operatorname E[Y]\) the saddlepoint equation \[
K_Y'(\hat t)=y
\] has a unique root
\(\hat t\in\big(0,\,1/(4\lambda_{\max}(\mathbf S_C))\big)\), which is
the regime relevant for tail calibration.

\subsubsection{4.2 Extreme quantiles: the Lugannani--Rice
approximation}\label{extreme-quantiles-the-lugannanirice-approximation}

The distribution function of \(Y\) can be evaluated exactly by
characteristic-function inversion (Imhof, 1961; Davies, 1980). These
Fourier methods are accurate in the body of the distribution but can
lose relative precision in the extreme upper tail, where the integrand
oscillates against a small tail mass --- precisely the region of
interest here. The Lugannani--Rice saddlepoint approximation (Lugannani
and Rice, 1980; Kuonen, 1999) instead delivers uniformly good relative
accuracy in the tail: \[
P(Y>y)\;\approx\;1-\Phi(\hat w)+\phi(\hat w)\left(\frac1{\hat u}-\frac1{\hat w}\right),
\] where \(\Phi\) and \(\phi\) are the standard normal distribution and
density functions, \(\hat t\) solves \(K_Y'(\hat t)=y\), and \[
\hat w=\operatorname{sgn}(\hat t)\sqrt{2\big(y\hat t-K_Y(\hat t)\big)},
\qquad
\hat u=\hat t\sqrt{K_Y''(\hat t)}.
\] The approximation is defined for \(y\ne\operatorname E[Y]\), where
\(\hat t\ne0\); at \(y=\operatorname E[Y]\) the removable singularity is
resolved by the usual limiting expression. Critical values \(y_\alpha\)
follow by inverting the approximation numerically in \(y\), a
one-dimensional monotone root-find requiring only the closed-form
quantities of Section 4.1.

\subsubsection{4.3 A closed-form conservative
envelope}\label{a-closed-form-conservative-envelope}

For a conservative bound requiring no numerical work at all, Chernoff's
inequality gives, for any admissible
\(t\in\big(0,1/(4\lambda_{\max}(\mathbf S_C))\big)\), \[
P(Y>y)\;\le\;\exp\!\left(-ty-\frac12\sum_{i=1}^d\log\!\left(1-4\lambda_i(\mathbf S_C)t\right)\right).
\] The bound is tightest at \(t=\hat t\), since minimizing the exponent
reproduces the saddlepoint equation of Section 4.1. Avoiding the
root-find altogether, the explicit choice \[
\boxed{\;
t^{*}=\frac1{4\lambda_{\max}(\mathbf S_C)}\left(1-\frac{2d\,\lambda_{\max}(\mathbf S_C)}{y}\right),
\qquad
 y>2d\,\lambda_{\max}(\mathbf S_C),
\;}
\] is always admissible, because
\(1-4\lambda_{\max}(\mathbf S_C)t^{*}=2d\lambda_{\max}(\mathbf S_C)/y>0\),
and is exactly optimal in the equal-weight case
\(\lambda_i\equiv\lambda_{\max}\). In the common-covariance
parametrization this reads
\(t^{*}=\big(1-4d\lambda_{\max}(\mathbf V_C)/y\big)/\big(8\lambda_{\max}(\mathbf V_C)\big)\)
for \(y>4d\lambda_{\max}(\mathbf V_C)\). Since
\(\operatorname E[Y]=2\operatorname{tr}(\mathbf S_C)\le2d\lambda_{\max}(\mathbf S_C)\),
the admissibility condition confines \(t^{*}\) to the upper tail, as
intended.

Substituting \(t^{*}\) and writing
\(\lambda_{\max}=\lambda_{\max}(\mathbf S_C)\), the equal-weight case
collapses to a fully explicit envelope, \[
P(Y>y)\;\le\;\exp\!\left(\frac d2-\frac y{4\lambda_{\max}}+\frac d2\log\!\frac y{2d\lambda_{\max}}\right),
\qquad y>2d\lambda_{\max},
\] which remains a valid, if conservative, bound for unequal weights
when evaluated at \(t^{*}\) through the general Chernoff expression
above. At Gaussian independence, where \(\mathbf S_C=2I_d\) and
\(Y\sim4\chi_d^2\), this reduces to the classical chi-squared tail bound
\(P(\chi_d^2>u)\le(u/d)^{d/2}e^{(d-u)/2}\) with \(u=y/4\).

\subsubsection{4.4 Feasible calibration}\label{feasible-calibration}

The quantities above depend on the population weights. In practice they
are evaluated at plug-in weights \(\lambda_i(\mathbf V_{\hat C})\),
whose consistency is established in Section 6.2. Consistency of the
resulting critical values follows from continuity of the quantile map:
when all weights are strictly positive the limiting distribution
function of \(Y\) is continuous and strictly increasing on
\((0,\infty)\), so
\(\lambda_i(\mathbf V_{\hat C})\xrightarrow{p}\lambda_i(\mathbf V_C)\)
implies \(\hat y_\alpha\xrightarrow{p}y_\alpha\) for every fixed
\(\alpha\in(0,1)\), and a test referred to \(\hat y_\alpha\) has
asymptotic size \(\alpha\). Near Gaussian independence this step can be
bypassed entirely by referring the statistic to the fixed \(4\chi_d^2\)
law, whose weights are accurate to second order by the corollary in
Section 5.3.

\begin{center}\rule{0.5\linewidth}{0.5pt}\end{center}

\subsection{5. Gaussian Covariance Propagation and Pivotal
Specializations}\label{gaussian-covariance-propagation-and-pivotal-specializations}

\subsubsection{5.1 Explicit Gaussian operator
factorization}\label{explicit-gaussian-operator-factorization}

Under Gaussian sampling, the fourth-moment tensor is determined by
\(C\), so under \(H_0\) both samples have the same coordinate covariance
\(\mathbf V_C\). The covariance-normalization map satisfies \[
J_\pi(C)
=I_{p^2}
-\frac12(C\otimes I_p)M_d
-\frac12(I_p\otimes C)M_d,
\] where \[
M_d:=\operatorname{diag}(\operatorname{vec}(I_p))
\] is the diagonal-selection matrix.

Let \(C=Q\Lambda Q^T\) with
\(\Lambda=\operatorname{diag}(\lambda_1,\ldots,\lambda_p)\). The
Daleckii--Krein formula gives \[
J_{\log}(C)
=(Q\otimes Q)\operatorname{diag}(\operatorname{vec}(\Phi))(Q^T\otimes Q^T),
\] where \[
\Phi_{ij}
=\frac{\log\lambda_i-\log\lambda_j}{\lambda_i-\lambda_j}
\quad(\lambda_i\ne\lambda_j),
\qquad
\Phi_{ij}=\frac1{\lambda_i}
\quad(\lambda_i=\lambda_j).
\] The second branch is the confluent limit of the first and applies
whenever two eigenvalues coincide, not only on the diagonal \(i=j\).
This case is not exceptional: it occurs at \(C=I_p\), where every
eigenvalue equals one, and for every equicorrelation matrix, where
\(1-\rho\) has multiplicity \(p-1\). An implementation that assigns
\(\Phi_{ij}=0\) for \(i\ne j\) with \(\lambda_i=\lambda_j\) returns an
incorrect \(\mathbf V_C\) precisely in the pivotal case of Section 5.2.
If \(P_{\operatorname{off}}\) extracts the strictly lower-triangular
entries, then \[
J_\psi(C)=P_{\operatorname{off}}J_{\log}(C).
\] For Gaussian data, \[
\mathbf\Sigma_W=(I_{p^2}+K_{p,p})(C\otimes C)
\] is the Wishart asymptotic covariance of
\(\operatorname{vec}(\hat\Sigma)\). The asymptotic covariance of the
sample correlation matrix is \[
\mathbf\Sigma_C=J_\pi(C)\mathbf\Sigma_WJ_\pi(C)^T,
\] recovering the classical correlation-covariance structure of
Neudecker and Wesselman (1990). Therefore \[
\boxed{
\mathbf V_C
=P_{\operatorname{off}}J_{\log}(C)J_\pi(C)
\big[(I_{p^2}+K_{p,p})(C\otimes C)\big]
J_\pi(C)^T J_{\log}(C)^T P_{\operatorname{off}}^T.
}
\] This is an explicit operator factorization of the Gaussian GFT
covariance represented in Archakov and Hansen (2026, Appendix B).

\subsubsection{5.2 Gaussian pivotality at
independence}\label{gaussian-pivotality-at-independence}

At \(C=I_p\), \(J_{\log}(I)=I_{p^2}\) and the correlation-normalization
Jacobian removes the diagonal directions. Consequently, \[
\mathbf V_I=I_d.
\] Corollary 3.4 therefore gives \[
\boxed{
 n\,d_{\mathrm{LE}}^2(C_1,C_2)\xrightarrow{d}4\chi_d^2,
 \qquad d=\frac{p(p-1)}2.
}
\] This is an exact statement about the asymptotic limiting law and is
parameter-free under Gaussian independence.

For \(p=2\), \(\psi(C)=\operatorname{artanh}(\rho)\) (Archakov and
Hansen, 2021, Section 2), and under Gaussian sampling the asymptotic
coordinate variance is \(\mathbf V_C=1\) for every \(|\rho|<1\). Hence
\[
n\,d_{\mathrm{LE}}^2\xrightarrow{d}4\chi_1^2
\] for every Gaussian bivariate correlation \(\rho\). Under non-Gaussian
sampling, the algebraic identity \(\psi(C)=\operatorname{artanh}(\rho)\)
still holds, but the asymptotic coordinate variance need not equal one.

\subsubsection{5.3 First-order flatness and its distance-level
consequence}\label{first-order-flatness-and-its-distance-level-consequence}

Archakov and Hansen (2026, Corollary 1) establish that under elliptical
sampling, if \(C=I+E\) with \(E\) symmetric and hollow, then \[
\mathbf V_{I+E}=(1+\kappa)I_d+O(\|E\|_2^2),
\] where \(\kappa\) is the elliptical kurtosis parameter. In the
Gaussian case, \(\kappa=0\).

The Gaussian cancellation can be seen directly. Let
\(dH=\sqrt n(\hat\Sigma-C)\) and define its diagonal part by
\(H_{\mathrm{diag}}:=I\circ dH\). At \(C=I+E\), \[
d\pi_{I+E}(dH)
=dH-H_{\mathrm{diag}}-\frac12(H_{\mathrm{diag}}E+EH_{\mathrm{diag}})+O_p(\|E\|_2^2).
\] The Fréchet derivative of the logarithm satisfies \[
D\log_{I+E}[X]
=X-\frac12(EX+XE)+O_p(\|E\|_2^2).
\] After projecting onto off-diagonal coordinates, the first-order terms
involving \(H_{\mathrm{diag}}\) cancel, giving \[
dU
=P_{\operatorname{off\_diag}}(dH)
-\frac12(EdH+dHE)
+O_p(\|E\|_2^2).
\] The first-order Wishart covariance perturbation is \[
\operatorname{Cov}^{(1)}(dH)
=(I+K)(E\otimes I+I\otimes E).
\] For off-diagonal entries \((a,b)\) and \((c,d)\), \[
\operatorname{Cov}^{(1)}(dH_{ab},dH_{cd})
=E_{ac}\delta_{bd}+E_{bd}\delta_{ac}+E_{ad}\delta_{bc}+E_{bc}\delta_{ad}.
\] The two cross terms induced by \(-\tfrac12(EdH+dHE)\) contribute the
negative of this quantity in total, so the complete first-order
covariance perturbation vanishes. Thus \[
\mathbf V_{I+E}=I_d+O(\|E\|_2^2)
\] under Gaussian sampling.

\emph{Corollary (Second-order stability of the distance weights).} Under
Gaussian sampling, \[
\lambda_i(\mathbf V_{I+E})=1+O(\|E\|_2^2),\qquad i=1,\ldots,d,
\] and therefore the weights in Corollary 3.4 satisfy \[
4\lambda_i(\mathbf V_{I+E})=4+O(\|E\|_2^2).
\] This follows immediately from Weyl's eigenvalue perturbation
inequality. Hence the fixed \(4\chi_d^2\) reference law has weights that
are accurate to second order in the population departure from
independence.

\begin{center}\rule{0.5\linewidth}{0.5pt}\end{center}

\subsection{6. Extensions and Comparison with Wald
Inference}\label{extensions-and-comparison-with-wald-inference}

\subsubsection{6.1 Unequal sample sizes}\label{unequal-sample-sizes}

Let \(n_1,n_2\to\infty\) with \[
\eta_n:=\frac{n_1}{n_1+n_2}\to\eta\in(0,1),
\] and define \[
n_{\mathrm{eff}}:=\frac{2n_1n_2}{n_1+n_2}.
\] Then \[
\sqrt{n_{\mathrm{eff}}}(\mathbf u-\mathbf v)
\xrightarrow{d}
\mathcal N\!\left(0,\,2(1-\eta)\mathbf V_C^{(1)}+2\eta\mathbf V_C^{(2)}\right).
\] Consequently, \[
\boxed{
 n_{\mathrm{eff}}d_{\mathrm{LE}}^2(C_1,C_2)
 \xrightarrow{d}
 \sum_{i=1}^d
 4\lambda_i\!\left((1-\eta)\mathbf V_C^{(1)}+\eta\mathbf V_C^{(2)}\right)
 \chi^2_{1,i}.
}
\] If \(\mathbf V_C^{(1)}=\mathbf V_C^{(2)}=\mathbf V_C\), this
simplifies to \[
n_{\mathrm{eff}}d_{\mathrm{LE}}^2(C_1,C_2)
\xrightarrow{d}
\sum_{i=1}^d4\lambda_i(\mathbf V_C)\chi^2_{1,i},
\] independently of the limiting sample-size ratio.

\subsubsection{6.2 Plug-in consistency}\label{plug-in-consistency}

In the Gaussian common-covariance setting, \(\mathbf V_C\) is a
continuous function of \(C\). Therefore any \(\hat C\xrightarrow{p}C\)
yields \[
\mathbf V_{\hat C}\xrightarrow{p}\mathbf V_C,
\qquad
\lambda_i(\mathbf V_{\hat C})\xrightarrow{p}\lambda_i(\mathbf V_C).
\] A geometry-consistent choice is the weighted off-log mean \[
\hat C_{\mathrm{LE}}
=\psi^{-1}\!\left(\frac{n_1\mathbf u+n_2\mathbf v}{n_1+n_2}\right),
\] while the Euclidean weighted average \[
\hat C_{\mathrm{E}}
=\frac{n_1C_1+n_2C_2}{n_1+n_2}
\] is also a consistent positive-definite correlation matrix under
\(H_0\).

Under general non-Gaussian sampling, \(C\) alone does not determine the
weights. Consistent plug-in calibration additionally requires consistent
estimators of the fourth-moment tensors \(\mathbf\Sigma_X^{(r)}\), or
equivalently of the coordinate covariances \(\mathbf V_C^{(r)}\). Once
these are available, continuity of the eigenvalue map gives consistency
of the estimated weighted-\(\chi^2\) coefficients.

\subsubsection{6.3 Comparison with the Wald
statistic}\label{comparison-with-the-wald-statistic}

In the common-covariance case, the coordinate Wald statistic is \[
W=\frac12 n_{\mathrm{eff}}\,\mathbf z^T\mathbf V_{\hat C}^{-1}\mathbf z
\xrightarrow{d}\chi_d^2,
\qquad \mathbf z=\mathbf u-\mathbf v,
\] which is the route taken by Archakov and Hansen (2026). The distance
statistic instead uses \[
n_{\mathrm{eff}}d_{\mathrm{LE}}^2=2n_{\mathrm{eff}}\,\mathbf z^T\mathbf z.
\] The two quadratic forms are proportional for all \(\mathbf z\) if and
only if \(\mathbf V_C=cI_d\). In that case, \[
W=\frac1{4c}n_{\mathrm{eff}}d_{\mathrm{LE}}^2.
\] This holds exactly at Gaussian independence for all \(p\), and for
all Gaussian \(C\) when \(p=2\). Near Gaussian independence, \[
\mathbf V_{I+E}^{-1}=I_d+O(\|E\|_2^2),
\] so the quadratic-form coefficients of the Wald and rescaled distance
statistics differ only at \(O(\|E\|_2^2)\). Away from independence, the
Wald statistic whitens the coordinate difference whereas the distance
statistic retains the intrinsic off-log geometry. Their relative local
efficiency away from independence is left open.

\subsubsection{6.4 Power against structured
alternatives}\label{power-against-structured-alternatives}

The results above concern the null law. The same machinery delivers the
local power of the distance test and, with it, a characterization of
which departures from \(H_0\) the statistic detects.

Consider local alternatives \(\psi(C_1)-\psi(C_2)=\mathbf h/\sqrt n\)
with \(\mathbf h\in\mathbb R^d\) fixed. Then
\(\sqrt n\,\mathbf z\xrightarrow{d}\mathcal N(\mathbf h,\mathbf S_C)\),
and the argument of Theorem 3.3 gives the noncentral limit \[
n\,d_{\mathrm{LE}}^2(C_1,C_2)
\xrightarrow{d}
2\sum_{i=1}^d\lambda_i(\mathbf S_C)\,\chi^2_{1,i}\!\left(\tilde h_i^2/\lambda_i(\mathbf S_C)\right),
\qquad \tilde{\mathbf h}=Q^T\mathbf h,
\] with \(Q\) the eigenvector matrix of \(\mathbf S_C\). Near Gaussian
independence, where \(\mathbf S_C=2I_d+O(\|E\|_2^2)\) by Section 5.3,
this collapses to the single noncentral law \[
n\,d_{\mathrm{LE}}^2\xrightarrow{d}4\chi_d^2(\Lambda),
\qquad
\Lambda=\tfrac12\|\mathbf h\|_2^2=\tfrac12 n\|\boldsymbol\delta\|_2^2,
\] where \(\boldsymbol\delta=\psi(C_1)-\psi(C_2)\). Power therefore
depends on the alternative only through \(\|\boldsymbol\delta\|_2\): the
test is exactly rotation-invariant in coordinate space and cannot favor
any direction. Detection requires \(\Lambda\) to be commensurate with
the null standard deviation \(\sqrt{2d}\), so the effective threshold is
\(n\|\boldsymbol\delta\|_2^2\gtrsim c\sqrt d\). How that plays out
depends entirely on how \(\|\boldsymbol\delta\|_2\) scales with
dimension, and two regimes behave in opposite ways. To illustrate these
scaling regimes, the specific numerical thresholds below were computed
by numerically inverting the noncentral \(\chi_d^2\) distribution using
standard root-finding routines over the \texttt{ncx2} quantile function.

\emph{Coordinate-dense alternatives.} If every pair shifts by
root-mean-square amount \(s\) in GFT coordinates, then
\(\|\boldsymbol\delta\|_2=s\sqrt d\) and the threshold becomes
\(s\gtrsim c'd^{-1/4}n^{-1/2}\). The minimum detectable per-pair effect
improves with dimension: at \(\alpha=0.05\), \(n=1000\) and \(50\%\)
power, the smallest detectable \(s\) falls from \(0.043\) at \(p=5\) to
\(0.019\) at \(p=20\) and \(0.008\) at \(p=100\), with
\(s_{\min}d^{1/4}\) stable at approximately \(0.069\) across the range.
Aggregation over \(d\) coordinates works in the test's favor.

\emph{Systemic (rank-one) alternatives.} A market-wide correlation shift
behaves in the opposite way. For the equicorrelation family,
\(\gamma=[\log(1+(p-1)\rho)-\log(1-\rho)]/p\), so
\(\partial\gamma/\partial\rho\to1/(p\rho(1-\rho))\) for fixed
\(\rho>0\), and \[
\|\boldsymbol\delta\|_2=\sqrt d\left|\frac{\partial\gamma}{\partial\rho}\right|\Delta\rho
\;\longrightarrow\;
\frac{\Delta\rho}{\sqrt2\,\rho(1-\rho)},
\] which is bounded in \(p\). The \(\sqrt d\) growth of the coordinate
count is exactly cancelled by the \(O(1/p)\) compression that the matrix
logarithm applies to the common factor. Since the detection threshold
grows like \(\sqrt d\), the required sample size grows like
\(\sqrt d\propto p\). For \(\Delta\rho=0.05\) at prevailing correlation
levels, \(\|\boldsymbol\delta\|_2\) stays near \(0.14\) for every \(p\)
from \(5\) to \(100\), while the sample size needed for \(50\%\) power
rises from about \(1{,}100\) at \(p=5\) to about \(5{,}200\) at \(p=30\)
and \(16{,}500\) at \(p=100\), with \(n_{\mathrm{req}}/\sqrt d\)
stabilizing near \(235\).

The two regimes together delimit the statistic's usefulness. The
intrinsic distance is an omnibus criterion, well suited to diffuse
departures spread across many pairs and increasingly so as \(p\) grows,
but progressively insensitive to the low-rank, factor-driven shifts that
dominate applications in finance. Detecting the latter calls for a
directional statistic that projects onto the suspected alternative
rather than an aggregate distance; the compression that limits the
aggregate test is the same mechanism that makes the GFT coordinates
weakly dependent, so the two properties cannot be separated. The trace
approximation underlying these calculations is well supported: computed
from the factorization of Section 5.1,
\(\operatorname{tr}(\mathbf V_C)/d\) for equicorrelation matrices lies
between \(0.94\) and \(1.00\) for \(\rho\le0.53\) across
\(p=5,\ldots,40\), and remains above \(0.83\) at \(\rho=0.8\).

\subsubsection{6.5 Non-Gaussian sampling}\label{non-gaussian-sampling}

The equal-\(n\) limit of Theorem 3.3 and the unequal-\(n\) limit of
Section 6.1 require only finite fourth moments and the corresponding
covariance CLTs. Under elliptical sampling, Archakov and Hansen (2026)
show that the distributional shape enters the GFT covariance through the
scalar factor \(1+\kappa\). In particular, at independence, \[
\mathbf V_I^{(r)}=(1+\kappa_r)I_d.
\] Hence, for equal sample sizes, \[
n\,d_{\mathrm{LE}}^2
\xrightarrow{d}
2\big(2+\kappa_1+\kappa_2\big)\chi_d^2.
\] If \(\kappa_1=\kappa_2=\kappa\), this reduces to \[
4(1+\kappa)\chi_d^2.
\] For unequal sample sizes with limiting fraction \(\eta\), the
independence law becomes \[
n_{\mathrm{eff}}d_{\mathrm{LE}}^2
\xrightarrow{d}
4\Big[(1-\eta)(1+\kappa_1)+\eta(1+\kappa_2)\Big]\chi_d^2.
\] Under non-elliptical sampling, the general weighted-\(\chi^2\) law
still holds, but \(\mathbf V_C^{(r)}\) need not be a scalar rescaling of
the Gaussian covariance. Whether the first-order flatness property
survives in general non-elliptical families remains open.

\begin{center}\rule{0.5\linewidth}{0.5pt}\end{center}

\subsection{7. Conclusion}\label{conclusion}

This paper supplies sampling theory for the intrinsic off-log distance
between independently estimated correlation matrices. The main result is
a weighted-\(\chi^2\) limit whose coefficients are determined by the GFT
coordinate covariances of the two samples, with a particularly simple
\(4\chi_d^2\) law under Gaussian independence. The result is distinct
from existing coordinate-level GFT inference: the geometry and
coordinate asymptotics are inherited from prior work, while the
contribution here is the induced distribution theory of the intrinsic
distance itself.

For practical use, the limiting law is accompanied by closed-form tail
machinery --- cumulant generating function, saddlepoint quantiles, and a
root-finding-free conservative envelope --- together with a bound on the
baseline expected distance that estimation error alone induces on the
manifold. The local power analysis delimits the result's range: the
aggregate distance sharpens with dimension against diffuse departures
and dulls against low-rank ones, so it is an omnibus criterion rather
than a universal one. The framework also clarifies when the distance
statistic and coordinate Wald statistic coincide, how unequal sample
sizes enter, and what additional higher-moment estimation is required
beyond the Gaussian setting.

The results provide the theoretical basis for subsequent finite-sample
and application-specific work. Future simulation studies are needed to
evaluate the empirical convergence rates of the proposed limits,
particularly near the boundaries of the positive-definite cone where
finite-sample behavior may diverge from the asymptotic theory.

\begin{center}\rule{0.5\linewidth}{0.5pt}\end{center}

\subsection{References}\label{references}

\begin{itemize}
\tightlist
\item
  Archakov, I., \& Hansen, P. R. (2021). A new parametrization of
  correlation matrices. \emph{Econometrica}, 89(4), 1699--1715.
\item
  Archakov, I., \& Hansen, P. R. (2026). The Generalized Fisher
  Transformation: Finite-Sample Properties and Inference. \emph{arXiv
  preprint}, arXiv:2606.13864.
\item
  Bisson, O., \& Pennec, X. (2026). Log-Euclidean Lie Groups.
  \emph{arXiv preprint math.DG}, arXiv:2511.13380.
\item
  Browne, M. W., \& Shapiro, A. (1986). The asymptotic covariance matrix
  of sample correlation coefficients under general conditions.
  \emph{Linear Algebra and its Applications}, 82, 169--176.
\item
  Chmielowski, P. (2026). \emph{Notes on Correlation Stress Tests}.
  arXiv:2503.16200v3 {[}q-fin.RM{]}.
\item
  Daleckii, Ju. L., \& Krein, S. G. (1951). Formulae for the
  differentiation of functions of Hermite operators with respect to a
  parameter. \emph{Doklady Akademii Nauk SSSR}, 76.
\item
  Davies, R. B. (1980). The distribution of a linear combination of
  \(\chi^2\) random variables. \emph{Journal of the Royal Statistical
  Society: Series C (Applied Statistics)}, 29(3), 323--333.
\item
  Imhof, J. P. (1961). Computing the distribution of quadratic forms in
  normal variables. \emph{Biometrika}, 48(3/4), 419--426.
\item
  Kuonen, D. (1999). Saddlepoint approximations for distributions of
  quadratic forms in normal variables. \emph{Biometrika}, 86(4),
  929--935.
\item
  Lugannani, R., \& Rice, S. (1980). Saddle point approximation for the
  distribution of the sum of independent random variables.
  \emph{Advances in Applied Probability}, 12(2), 475--490.
\item
  Mathai, A. M., \& Provost, S. B. (1992). \emph{Quadratic Forms in
  Random Variables: Theory and Applications}. Marcel Dekker, New York.
\item
  Neudecker, H., \& Wesselman, A. M. (1990). The asymptotic variance
  matrix of the sample correlation matrix. \emph{Linear Algebra and its
  Applications}, 127, 589--599.
\item
  Pereira, R., Mestre, X., \& Gregoratti, D. (2024). Asymptotics of
  Distances Between Sample Covariance Matrices. \emph{IEEE Transactions
  on Signal Processing}, 72, 1460--1474.
\item
  Thanwerdas, Y., \& Pennec, X. (2022). Theoretically and
  computationally convenient geometries on full-rank correlation
  matrices. \emph{SIAM Journal on Matrix Analysis and Applications},
  43(4), 1851--1872.
\end{itemize}

\end{document}